\documentclass[10pt,reqno]{amsart}
\usepackage{graphicx}
\usepackage{epic}
\usepackage{eepic}
\usepackage{url}

\theoremstyle{plain}
   \newtheorem{theorem}{Theorem}[section]
   \newtheorem{proposition}[theorem]{Proposition}
   
   \newtheorem{corollary}[theorem]{Corollary}
   \newtheorem{conjecture}[theorem]{Conjecture}
   
\theoremstyle{definition}
   \newtheorem{definition}[theorem]{Definition}

   \newtheorem{remark}[theorem]{Remark}

\author[P.~Br\"and\'en]{Petter Br\"and\'en}
       \address{Department of Mathematics, University of Michigan,
       Ann Arbor, MI 48109-1043, USA}
       \email{branden@umich.edu}
       

\numberwithin{equation}{section}
\def\dd{\kern.4ex\mbox{\raise.4ex\hbox{{\rule{.35em}{.12ex}}}}\kern.4ex}
\def\132{\ensuremath{(13{\dd}2)}}
\def\231{\ensuremath{(2{\dd}31)}}
\newcommand{\NN}{\mathbb{N}}
\newcommand{\ZZ}{\mathbb{Z}}

\newcommand{\sym}{\mathfrak{S}}
\newcommand{\JH}{\mathcal{L}}
\newcommand{\WP}{\overline{W}}

\def\newop#1{\expandafter\def\csname #1\endcsname{\mathop{\rm
#1}\nolimits}}

\newop{des}
\newop{peak}
\newop{v}
\newop{cpeak}
\newop{veh}
\newop{SIVEH}
\newop{Redge}
\newop{Odd}
\newop{EV}
\newop{MAJ}
\newop{Orb}
\newop{Odd}
\newop{Des}
\newcommand{\vej}{\veh'}

\begin{document}


%

%

\title[Actions on Permutations]
{Actions on Permutations and Unimodality of 
Descent Polynomials}
\begin{abstract}
We study a group action on permutations due to Foata and Strehl and use it to prove that the descent generating polynomial of certain sets of permutations has a nonnegative 
expansion in the basis $\{ t^i(1+t)^{n-1-2i}\}_{i=0}^m$, $m=\lfloor (n-1)/2 \rfloor$. This property implies symmetry and unimodality. We prove that the action is invariant under stack-sorting which strengthens recent unimodality results 
of B\'ona. We prove that the generalized permutation patterns $\132$ and $\231$ are invariant under the action and use this to prove unimodality properties for a $q$-analog of the Eulerian numbers 
recently studied by Corteel, Postnikov, Steingr\'{\i}msson and Williams.   

We also extend the action to linear extensions of sign-graded posets to give a new proof of the unimodality of the $(P,\omega)$-Eulerian polynomials of sign-graded posets and a combinatorial 
interpretations (in terms of Stembridge's peak polynomials) of the corresponding coefficients when expanded in the above basis.  

Finally, we prove that the statistic defined as the number of vertices of even height in the unordered decreasing tree of a permutation has the same distribution as the number of descents on any 
set of permutations invariant under the action. When restricted to the set of  stack-sortable permutations we recover a result of Kreweras. 
\end{abstract}
\maketitle
\thispagestyle{empty}
\section{Introduction}
The $n$-th {\em Eulerian polynomial}, $A_n(t)= A_{n1} + A_{n2}t+\cdots+A_{n(n-1)}t^{n-1}$, may be defined as the generating polynomial for 
the number of {\em descents} over the symmetric group $\sym_n$, i.e., 
$$
A_n(t)=\sum_{\pi \in \sym_n}t^{\des(\pi)}, 
$$
where $\des(\pi)=|\{ i : a_i > a_{i+1}\}|$ and where $\pi : i \rightarrow a_i$ ($1 \leq i \leq n$) is identified 
with the word $a_1a_2\cdots a_n$ in the distinct $n$ letters $a_1, \ldots, a_n$ taken out of 
$[n]:=\{1,2,\ldots, n\}$.  

 

In a series of papers \cite{Foata, Foata-Strehl-1,Foata-Strehl-2} Foata and Strehl studied a group action on the symmetric group, $\sym_n$, with the following properties. 
The number of orbits is the $n$-th tangent number or secant number, 
according as $n$ is odd or even, and if an orbit, $\Orb(\pi)$, of a permutation $\pi \in \sym_n$ is enumerated according to the number of descents then 
\begin{equation}\label{fs}
\sum_{\sigma \in \Orb(\pi)}t^{\des(\sigma)}= (2t)^{\v(\pi)}(1+t)^{n-1-2\v(\pi)}, 
\end{equation}
where $\v(\pi)=|\{i: a_{i-1}>a_i<a_{i+1}\}|$. 
 From \eqref{fs} it follows that $A_n(t)$ has nonnegative coefficients when expanded in the basis $\{ t^k(1+t)^{n-1-2k}\}_{k=0}^{\lfloor (n-1)/2 \rfloor}$,  a
result which can also be proven analytically \cite{Carlitz-Scoville,Foata-Schutzenberger}. This implies 
that the sequence  $\{A_{ni}\}_{i=0}^{n-1}$ is symmetric and unimodal, i.e.,  that $A_{ni}=A_{n(n-1-i)}$, $1\leq i \leq n-1$ and
$$
A_{n0} \leq A_{n1} \leq \cdots \leq A_{nc} \geq A_{n(c+1)} \geq \cdots \geq A_{n(n-1)},
$$
where $c=\lfloor (n-1)/2 \rfloor$. Indeed, $\{A_{ni}\}_{i=0}^{n-1}$ is a nonnegative sum of unimodal and 
symmetric sequences with the same center of symmetry. 

We will in this paper study a slightly modified version of the Foata-Strehl action
and show that interesting subsets of $\sym_n$ are invariant under the action. In particular we show that the set of $r$-stack sortable permutations is invariant under the action which strengthens the recent result of B\'ona \cite{Bona-Uni,Bona-Cor} claiming that the corresponding descent generating polynomial  is symmetric and unimodal.  

 In Section~\ref{pq} we prove that the generalized permutation patterns $\132$ and $\231$ are invariant under the modified Foata-Strehl action. This is used to prove unimodality properties for a $q$-analog of the Eulerian numbers 
recently studied by Corteel, Postnikov, Steingr\'{\i}msson and Williams \cite{Corteel,Cort-Will,Postnikov,Stein-Will,Williams} and which appears as a translation of the polynomial enumerating  the cells in the totally nonnegative part of a Grassmannian  
\cite{Postnikov,Williams}, and also in the stationary distribution of the ASEP model in statistical mechanics \cite{Corteel,Cort-Will}.

We will in Section~\ref{sign-graded} define an action on the set of linear extensions of  a sign-graded poset, see Section~\ref{sign-graded} for relevant definitions. This enables us  to give a combinatorial interpretation in terms of Stembridge's peak polynomials of the coefficients of the $(P,\omega)$-Eulerian polynomials when expanded in the basis 
$\{t^i(1+t)^{d-2i}\}_{i=0}^{\lfloor d/2 \rfloor}$, $d=|P|-r-1$. 

%

In Section~\ref{veh} we study the statistic $\pi \rightarrow \veh(\pi)$ on permutations which is defined as the number of vertices of even height in the unordered increasing tree of $\pi$. We prove that $\veh$ has the same distribution as $\des$ on every subset of $\sym_n$ invariant under the action.  This 
can be seen as a generalization of a result of Kreweras \cite{Kreweras}. In Section~\ref{mahonian} we also find a Mahonian partner for $\veh$. 

Finally, in Section~\ref{OP}, we discuss further directions and open problems. 
\section{The Action of Foata and Strehl}
Let $\pi=a_1a_2\cdots a_n \in \sym_n$ and let $x  \in [n]$. We may write $\pi$ as the concatenation 
$\pi = w_1w_2xw_4w_5$ where $w_2$ is the maximal contiguous subword immediately to the left 
$x$ whose letters are all smaller than $x$, and $w_4$ is the maximal contiguous subword immediately to the right of $x$ whose letters are all smaller than $x$. This is the $x$-{\em factorization} of $\pi$. Define $\varphi_x(\pi) = w_1w_4xw_2w_5$. Then $\varphi_x$ is an involution acting on $\sym_n$ and 
it is not hard to see that $\varphi_x$ and $\varphi_y$ commute for all $x,y \in [n]$. Hence 
for any subset $S \subseteq [n]$ we may define the 
function $\varphi_S : \sym_n \rightarrow \sym_n$ by
$$
\varphi_S(\pi) = \prod_{x \in S} \varphi_x(\pi).
$$
The group $\ZZ_2^n$ acts on $\sym_n$ via the 
functions $\varphi_S$, $S \subseteq [n]$. This action was studied by Foata and Strehl in 
\cite{Foata, Foata-Strehl-1,Foata-Strehl-2}. To be precise, Foata and Strehl defined the action as 
$C \circ \varphi_S \circ C$, where $C : \sym_n \rightarrow \sym_n$ is the involution described by 
$a_1a_2\cdots a_n \mapsto b_1b_2\cdots b_n$, where $b_i= n+1-a_i$, $1\leq i \leq n$.  Sometimes it is preferable to define the action on the {\em decreasing binary tree} of the permutation. The decreasing binary tree of a permutation of a finite subset of $\{1,2,3,\ldots\}$ is defined recursively as follows. The empty tree corresponds to the empty word. If 
$\pi$ is non-empty then we may write $\pi$ as the concatenation $\pi=LmR$ where $m$ and $L$ and $R$ are the subwords to the left and right of $m$ respectively. The tree corresponding to $\pi$ has a root labeled $m$ and as left subtree the tree corresponding to $L$ and as right subtree the tree corresponding to $R$. This describes a bijective correspondence between the set of decreasing binary trees with labels $[n]$ and $\sym_n$. It is not hard to see that the tree of $\varphi_x(\pi)$ is obtained by exchanging the subtrees rooted at $x$, if any. Another action on permutations with similar properties was studied by Hetyei and Reiner \cite{HR} and subsequently by Foata and Han \cite{FH}.

\begin{figure}\caption{\label{graph} Graphical representation of $\pi = 573148926$. The dotted lines indicates where 
the double ascents/descents move to.}
\bigskip
\setlength{\unitlength}{5mm}
\begin{picture}(9,11)
\put(0,0){
\path(0,10)(1,5)(2,7)(3,3)(4,1)(5,4)(6,8)(7,9)(8,2)(9,6)(10,10)
\put(0,10){$10$}\put(1,4.5){5}\put(2,7.1){7}\put(3.1,3){3}\put(4,.5){1}
\put(5.2,4){4}\put(5.7,8){8}\put(7,9){9}\put(8,1.5){2}\put(9.2,6){6}
\put(10,10){$10$}
\dottedline{0.3}(7.5,6)(9,6)
\dottedline{0.3}(3,4)(5,4)
\dottedline{0.3}(3,3)(4.5,3)
\dottedline{0.3}(0.7,8)(6,8)
}
\end{picture}
\end{figure}
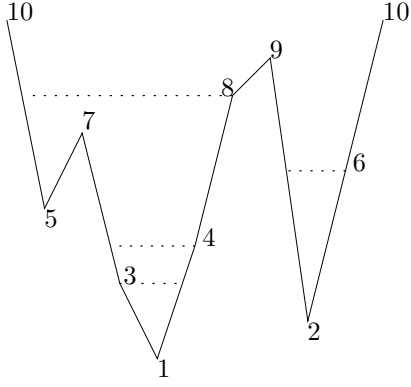

Let $\pi = a_1a_2\cdots a_n$ be a permutation in $\sym_n$ 
and let $a_0=a_{n+1}=n+1$. 
If $k \in [n]$ then $a_k$ is a
\begin{itemize}
\item[] {\em valley} if $a_{k-1}> a_k < a_{k+1}$,
\item[] {\em peak} if $a_{k-1}<  a_k >a_{k+1} $,
\item[] {\em double ascent} if  $a_{k-1} <a_k < a_{k+1}$, and 
\item[] {\em double descent} if $a_{k-1}> a_k >a_{k+1}$.
\end{itemize}
Let $x \in [n]$ and let $\pi=a_1a_2\ldots a_n \in \sym_n$. We make the following observation. 
\begin{itemize}
\item If $x$ is a double descent then $\varphi_x(\pi)$ is 
obtained by inserting $x$ between the first pair of letters 
$a_i, a_{i+1}$ to the right of $x$ such that $a_i < x < a_{i+1}$.
 
\item If $x$ is a double ascent then $\varphi_x(\pi)$ is 
obtained by inserting $x$ between the first pair of letters 
$a_i, a_{i+1}$ to the left of $x$ such that $a_i > x > a_{i+1}$.
\end{itemize}
We modify the Foata-Strehl action in the following way. If $x \in [n]$ then 
$$
\varphi'_x(\pi)= 
\begin{cases} 
\varphi_x(\pi) \mbox{ if } x \mbox{ is a double ascent or double descent,} \\
\pi \ \ \ \ \ \ \mbox{ if } x \mbox{ is a valley or a peak.}
\end{cases}
$$
The functions are easily visualized when a permutation is represented 
graphically. Let $\pi=a_1a_2\cdots a_n \in \sym_n$ and imagine 
marbles at coordinates $(i,a_i)$, $i=0,1,\ldots,n+1$ in the grid $\NN \times \NN$. For $i = 0,1,\ldots, n$ 
connect $(i,a_i)$ and $(i+1,a_{i+1})$ with a wire. Suppose that 
gravity acts on the marbles from above and  suppose that $x$ is not at an equilibrium. If we release $x$ from the wire it will slide and stop when 
it has reached the same height again. The resulting permutation will be  
$\varphi'_x(\pi)$, see Fig.~\ref{graph}.  The functions $\varphi'_x$ were studied by Shapiro, Woan and Getu unaware\footnote{The present author was also unaware of this until it was pointed out by the referee.} that they are essentially the same as the functions defining the Foata-Strehl action.

Again it is clear that the $\varphi'_x$'s are involutions and that they commute. Hence, 
for any subset $S \subseteq [n]$ we may define the 
function $\varphi'_S : \sym_n \rightarrow \sym_n$ by
$$
\varphi'_S(\pi) = \prod_{x \in S} \varphi'_x(\pi).
$$
Hence the group $\ZZ_2^n$ acts on $\sym_n$ via the 
functions $\varphi'_S$, $S \subseteq [n]$. Subsequently we will refer to this action as the modified Foata-Strehl action, or the MFS-action for short. 
\section{Properties of the Modified Foata-Strehl Action} 
For $\pi \in \sym_n$ let $\Orb(\pi)= \{ g(\pi) : g \in \ZZ_2^n\}$ be the 
orbit of $\pi$ under the MFS-action. There is a unique element in $\Orb(\pi)$ 
which has no double descents and which we denote by  $\hat{\pi}$. The next theorem follows from the work  in
\cite{Foata-Strehl-2,SWG}, but we prove it here for completeness.
\begin{theorem}\label{orb}
Let $\pi \in \sym_n$. Then 
$$
\sum_{\sigma \in \Orb(\pi)} t^{\des (\sigma) } = 
t^{\des (\hat{\pi})} (1+t)^{n-1-2\des (\hat{\pi})}=t^{\peak (\pi)} (1+t)^{n-1-2\peak(\pi)},
$$
where $\peak(\pi)=|\{ i : a_{i-1}<a_i>a_{i+1}\}|$.
\end{theorem}
\begin{proof}
If $x$ is a double ascent in $\pi$ then 
$\des (\varphi'_x(\pi))= \des (\pi) +1$. It follows that 
$$
\sum_{\sigma \in \Orb(\pi)} t^{\des (\sigma) }= t^{\des(\hat{\pi})}(1+t)^a,
$$
where $a$ is the number of double ascents in $\hat{\pi}$. If we 
delete all double descents from $\hat{\pi}$ we get 
an alternating permutation 
$$
n+1 > b_1 < b_2 > b_3 <\cdots > b_{n-a}< n+1,
$$
with the same number of descents. Hence  
$n-a = 2\des (\hat{\pi}) +1$. Clearly $\des(\hat{\pi})=\peak(\pi)$ and the theorem follows. 
\end{proof}  
For a subset $T$ of $\sym_n$ let 
$$
W(T;t)=\sum_{\pi \in T} t^{\des(\pi)}   \ \ \mbox{ and } \ \   \WP(T;t)=\sum_{\pi \in T} t^{\peak(\pi)}.
$$
\begin{corollary}\label{corre} 
Suppose that $T \subseteq \sym_n$ is invariant under the MFS-action. Then 
$$
W(T;t) = 2^{-n+1}(1+t)^{n-1}\WP(T; 4t(1+t)^{-2}).
$$
Equivalently
$$
W(T;t) = \sum_{i=0}^{\lfloor n/2 \rfloor} b_i(T)t^i(1+t)^{n-1-2i}, 
$$
where 
$$
b_i(T)= 2^{-n+1+2i}|\{ \pi \in T : \peak(\pi)=i\}|
$$
\end{corollary}
\begin{proof}
It is enough to prove the theorem for an orbit of a permutation $\pi \in \sym_n$. Since the number of  peaks is constant on $\Orb(\pi)$ the equality follows from Theorem~\ref{orb}. 
\end{proof}
\begin{remark}\label{comb}
If  we want to prove "combinatorially" that the coefficients of $W(T;t)$ form a symmetric and unimodal sequence then we can construct an involution  proving symmetry and an injection proving unimodality easily as follows. 

Define  $f: \sym_n \rightarrow \sym_n$ by $f =  \varphi'_{[n]}$. Clearly $f$ is an involution and restricts to any subset of $\sym_n$ invariant under the MFS-action. 
Moreover,
\begin{equation}\label{symf}
\des(f(\pi))+ \des(\pi)=n-1, 
\end{equation}
so $f$ has the desired properties. The involution $f$ was defined differently in \cite{Bona-Uni}. To find an injection 
$$
g_j : \{ \pi \in T : \des(\pi) =j \}  \rightarrow \{\pi \in T : \des(\pi)=j+1\}, 
$$
for $j =1,2,\ldots, \lfloor (n-1)/2 \rfloor$  it suffices to find an injection from the set of subsets of cardinality 
$k$ of $[m]$ to the set of subsets of cardinality $k+1$ of $[m]$, for $1 \leq k \leq \lfloor m/2\rfloor$. This can 
done as in e.g. \cite{Sagan}. 
\end{remark}

\section{Invariance Under Stack Sorting}
Much has been written on the combinatorics of the stack-sorting problem (cf.~\cite{Bona-Survey}) since it was introduced by Knuth \cite{Knuth}.  
The stack-sorting operator $S$ can be defined recursively on permutations of finite subsets of 
$\{1,2,\ldots \}$ as follows. If 
$w$ is empty then $S(w)=w$ and  if  
$w$ is non-empty 
write $w$ as the concatenation $w=LmR$, where $m$ is the greatest element 
of $w$ and $L$ and $R$ are the subwords to the left and right of $m$ 
respectively. Then $S(w)=S(L)S(R)m$. 

Let $\pi = a_1a_2\cdots a_n \in \sym_n$. Recall that $i \in [n-1]$ is a descent in $\pi$ if $a_i>a_{i+1}$. If $i$ is a descent in $\pi$ we 
let $r_i(\pi)$ be the permutation 
obtained by inserting $a_i$ between the first pair of letters 
$a_j, a_{j+1}$ to the right of $x$ such that $a_j < x < a_{j+1}$ ($a_{n+1}= n+1$). The following theorem 
describes a new way of computing $S(\pi)$.  

\begin{theorem}
Let $i_1<i_2 < \cdots <i_d$ be the descents in the permutation 
$\pi=a_1a_2\cdots a_n$. Then 
$$
S(\pi) = r_{i_d}r_{i_{d-1}}\cdots r_{i_1}(\pi).
$$
\end{theorem}
\begin{proof}
Let $S' : \sym_n \rightarrow \sym_n$ be defined by $S'(\pi) = r_{i_d}r_{i_{d-1}}\cdots r_{i_1}(\pi)$. 
It is straightforward to check that $S'$ satisfies the same recursion as $S$. 
\end{proof}
\begin{figure}\caption{\label{stacksort} Computing $S(573148926)=r_7r_3r_2(573148926)=513478269$. }
\bigskip
\setlength{\unitlength}{4mm}
\begin{picture}(9,11)
\put(0,0){
\path(0,10)(1,5)(2,7)(3,3)(4,1)(5,4)(6,8)(7,9)(8,2)(9,6)(10,10)
\put(0,10){$10$}\put(1,4.5){5}\put(2,7.1){7}\put(3.1,3){3}\put(4,.5){1}
\put(5.2,4){4}\put(5.7,8){8}\put(7,9){9}\put(8,1.5){2}\put(9.2,6){6}
\put(10,10){$10$}
\put(11.5,6.7){$r_2$}
\put(11,6){$\longrightarrow$}
\dottedline{0.3}(2,7)(5.8,7)

}
\end{picture}
\ \ \ \ \ \ \ \ \ \ \ \ 
\begin{picture}(9,11)
\put(0,0){
\path(0,10)(1,5)(2,3)(3,1)(4,4)(5,7)(6,8)(7,9)(8,2)(9,6)(10,10)
\put(0,10){$10$}\put(1,4.5){5}\put(4.8,7){7}\put(2.1,3){3}\put(3,.5){1}
\put(4.2,4){4}\put(5.7,8){8}\put(7,9){9}\put(8,1.5){2}\put(9.2,6){6}
\put(10,10){$10$}
\put(11.5,6.7){$r_3$}
\put(11,6){$\longrightarrow$}
\dottedline{0.3}(2,3)(3.8,3)
}
\end{picture}
\ \ \ \ \ \ \ \ \ \ \ \ 
\begin{picture}(9,11)
\put(0,0){
\path(0,10)(1,5)(2,1)(3,3)(4,4)(5,7)(6,8)(7,9)(8,2)(9,6)(10,10)
\put(0,10){$10$}\put(1,4.5){5}\put(4.8,7){7}\put(3.1,3){3}\put(2,.5){1}
\put(4.2,4){4}\put(5.7,8){8}\put(7,9){9}\put(8,1.5){2}\put(9.2,6){6}
\put(10,10){$10$}
\put(11.5,6.7){$r_7$}
\put(11,6){$\longrightarrow$}
\dottedline{0.3}(7,9)(9.8,9)
}
\end{picture}
\ \ \ \ \ \ \ \ \ \ \ \ 
\begin{picture}(9,11)
\put(0,0){
\path(0,10)(1,5)(2,1)(3,3)(4,4)(5,7)(6,8)(7,2)(8,6)(9,9)(10,10)
\put(0,10){$10$}\put(1,4.5){5}\put(4.8,7){7}\put(3.1,3){3}\put(2,.5){1}
\put(4.2,4){4}\put(5.7,8){8}\put(9,9){9}\put(7,1.5){2}\put(8.2,6){6}
\put(10,10){$10$}
}
\end{picture}
\end{figure}
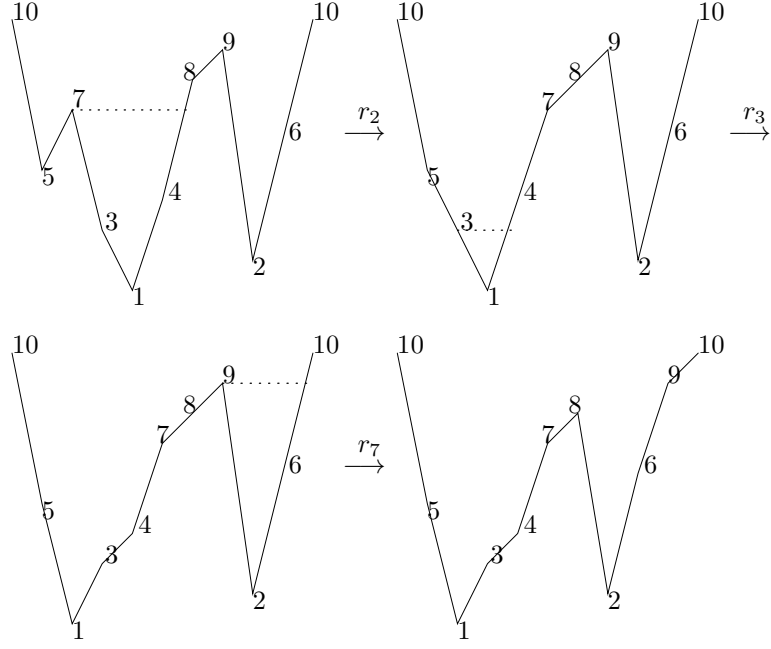
From the above description of $S$ we see that $S( \varphi_x (\pi) )= S(\pi)$, hence the following corollary. 
\begin{corollary}\label{stackop}
If $\sigma, \tau \in \sym_n$ are in the same orbit under 
the MFS-action then $S(\sigma)=S(\tau)$. 
\end{corollary}
Corollary~\ref{stackop} can also be deduced from \cite[Proposition 2.1]{Bousquet}.

Let $r \in \NN$. 
A permutation $\pi \in \sym_n$ is said to be $r${-stack sortable} if  
$S^r(\pi)=12\cdots n$. Denote by $\sym_n^r$ the set of $r$-stack sortable permutations in $\sym_n$. By Corollary~\ref{stackop} we have that $\sym_n^r$ is invariant under the MFS-action for all $n,r \in \NN$ so Corollary~\ref{corre} applies.
\begin{corollary}\label{genbona}
For all $n,r \in \NN$ we have 
$$
W(\sym_n^r;t) = \sum_{i=0}^{\lfloor n/2 \rfloor} b_i(\sym_n^r)t^i(1+t)^{n-1-2i}, 
$$
where 
$$
b_i(\sym_n^r)= 2^{-n+1+2i} |\{ \pi \in \sym_n^r : \peak(\pi)=i\}|.
$$
\end{corollary}
An immediate consequence of Corollary~\ref{genbona} is the following theorem due to 
B\'ona. 
\begin{theorem}[B\'ona \cite{Bona-Uni,Bona-Cor}] 
For all $n,r \in \NN$, the coefficients of $W(\sym_n^r;t)$ form a  symmetric and unimodal sequence.
\end{theorem}
An open problem posed by B\'ona 
\cite{Bona-Uni} is to determine whether the polynomial $W(\sym_n^r;t)$ has the stronger property of having all zeros real for $n,r \in \NN$. This is known for $r \geq n-1$ because then $W(\sym_n^r;t)=A_n(t)$ and the Eulerian polynomials are known to have all zeros real (cf.~\cite{Harper}), and for $r=1$ as we then get the Narayana polynomials \eqref{nara} 
which are known to have all zeros real by e.g. Malo's theorem (cf.~\cite{Marden}).  In 
\cite{Branden-Trans} we prove real-rootedness whenever $r=2$ or $r=n-2$. It is easy to see (cf.~\cite{Branden-Sign}) that 
if all zeros of $p(t)=\sum_{i=0}^n a_i t^i$ are real and $\{a_i\}_{i=0}^n$ is nonnegative and symmetric with center of symmetry $d/2$, then 
$$
p(t) = \sum_{i=0}^{\lfloor d/2 \rfloor} b_i t^i(1+t)^{d-2i},
$$ 
where $b_i$, $i=0,\ldots, \lfloor d/2 \rfloor$ are nonnegative. Hence Corollary~\ref{genbona} can be seen as further evidence for a positive answer to B\'ona's question.

Knuth  \cite{Knuth} proved that the $1$-stack sortable permutations are exactly the permutations that avoid the pattern 
$231$, i.e.,  permutations $\pi = a_1a_2\cdots a_n$ such that $a_k<a_j<a_i$ for no $1\leq i < j < k\leq n$. The set of  $231$-avoiding permutations in $\sym_n$ is denoted by $\sym_n(231)$. Simion \cite{Simion} proved that the $n$-th Narayana polynomial is the descent generating polynomial of $\sym_n(231)$, i.e.,    
\begin{equation}\label{nara}
\begin{split}
W(\sym_n(231);t)&=\sum_{k=0}^{n-1} \frac1 n \binom n k \binom n {k+1}t^k \\
&= \sum_{k=0}^{\lfloor n/2 \rfloor} \frac 1 {k+1} \binom {2k} k \binom {n-1}{2k} t^k(1+t)^{n-1-2k}, 
\end{split}
\end{equation}
where the second equality can be derived using hypergeometric formulas, see also \cite{Simion-Ullman}. Hence we have the following corollary. 
\begin{corollary}
Let $n,k \in \NN$. Then 
$$
|\{ \pi \in \sym_n(231): \peak(\pi)=k\}| = 2^{n-1-2k}\frac {1} {k+1} \binom {2k} k \binom {n-1}{2k}.
$$
\end{corollary}

\section{A Refinement of the Eulerian Polynomials}\label{pq}
The statistic $\231 : \sym_n \rightarrow \sym_n$ is an instance of a generalized permutation pattern as introduced by Babson and Steingr\'{\i}msson \cite{Babson}. Let 
$\pi =a_1a_2\cdots a_n \in \sym_n$. Then $\231(\pi)$ is the number  
of pairs  $1 \leq i < j \leq n-1$ such that $a_{j+1}<a_i<a_j$. Similarly,  
let $\132(\pi)$ be the number of pairs $2 \leq i < j \leq n$ such that 
$a_{i-1}<a_j<a_i$. 
\begin{theorem}\label{constant}
The statistics $\231$ and $\132$ are constant on any orbit under 
the MFS-action.
\end{theorem}
\begin{proof}
An alternative description of $\231(\pi)$, $\pi=a_1a_2\cdots a_n$ is the number triples 
$(a_i, a_j, a_k)$ such that $1\leq i < j<k \leq n$ and $a_k<a_i<a_j$, where $(a_j,a_k)$ is a pair of  consecutive 
peak and valley. By consecutive we mean that there are no other peaks or valleys in between $a_j$ and $a_k$. The number of such triples is invariant under the action since $a_j$ and $a_k$ cannot move and $a_i$ cannot move over the peak $a_j$. A similar reasoning applies to $\132$. 
\end{proof}

Define a $(p,q)$-refinement of the Eulerian polynomial by
$$
A_n(p,q,t)= \sum_{\pi \in \sym_n}p^{\132(\pi)}q^{\231(\pi)}t^{\des (\pi)}.
$$
These polynomials (or at least $A_n(p,1,t)$ and $A_n(1,q,t)$) have been in focus in several recent papers \cite{Corteel,Cort-Will,Postnikov,Stein-Will,Williams}. A fascinating property of the polynomial $A_n(p,1,t)$ is that it appears as a translation of the polynomial enumerating  the cells in the totally nonnegative part of a Grassmannian  
\cite{Postnikov,Williams}, and also in the stationary distribution of the ASEP model in statistical mechanics \cite{Corteel,Cort-Will}. 
 
From Theorem~\ref{constant} and Theorem~\ref{corre} we get that
$$
A_n(p,q,t)= \sum_{i=0}^{\lfloor (n-1)/2 \rfloor} 
b_{n,i}(p,q)t^i(1+t)^{n-1-2i},
$$
where 
\begin{equation}\label{bni}
b_{n,i}(p,q)= 2^{-n+1+2i}\mathop{\sum_{\pi \in \sym_n}}_{\peak(\pi)=i} p^{\132(\pi)}q^{\231(\pi)}.
\end{equation}
\begin{proposition}
Let $n \in \NN$. Then 
$$
A_n(p,q,t)=A_n(q,p,t).
$$
\end{proposition}

\begin{proof}
Let $f$ be as in Remark~\ref{comb} and let $R: \sym_n \rightarrow \sym_n$ be defined  by 
$$
R(\pi)= a_n \cdots a_2a_1, \ \ \mbox{ if }  \pi=a_1a_2 \cdots a_n. 
$$
Let $\pi'=R(f(\pi))$. Then 
$$
\Big(\des(\pi'), \132(\pi'), \231(\pi')\Big)= \Big(\des(\pi), \231(\pi),\132(\pi)\Big), 
$$
and the proposition follows. 
\end{proof}

A further striking property of $A_n(p,q,t)$ is that 
$$
A_n(q,q^2,q)=A_n(q^2,q,q)= [n]_q[n-1]_q \cdots [1]_q, 
$$
where $[k]_q= 1+q+q^2+\cdots+q^{k-1}$. 
This is because the statistics
$$
S_1= \132+\132 +\231 +\des \ \ \mbox{ and } \ \ S_2= \132+\231 +\231 +\des 
$$
are Mahonian (see Section~\ref{mahonian}), a fact due to Simion and Stanton \cite{Simion-Stanton}, see also \cite{Babson}.

\section{An Action on the Linear Extensions of a Sign-Graded Poset}\label{sign-graded}
Recall that a {\em labeled poset} is a pair $(P,\omega)$ where $P$ is a 
finite poset and $\omega : P \rightarrow \ZZ$ is an injection. The 
{\em Jordan-H\"older set}, $\JH(P,\omega)$, is the set of permutations 
$\pi = a_1a_2\cdots a_p$ ($p=|P|$) of $\omega(P)$ such that if $x$ is smaller than 
$y$ in $P$ ($x <_P y$), then 
$\omega(x)$ precedes $\omega(y)$ in $\pi$. The 
$(P,\omega)$-{\em Eulerian polynomial} is defined 
by 
$$
W(P,\omega;t)=\sum_{\pi \in \JH(P,\omega)}t^{\des (\pi)}. 
$$
Hence the $n$-th Eulerian polynomial is the $(P,\omega)$-Eulerian polynomial of an anti-chain of size $n$. 
The $(P,\omega)$-Eulerian polynomials have been intensively studied since they were introduced by 
Stanley \cite{Stanley-thesis} in 1972. For example, the Neggers-Stanley conjecture which asserts that these polynomials always 
have real zeros has attracted widespread attention \cite{Athanasiadis,Bjorner-Farley,Branden-Diamond,Branden-Sign,Branden-CE,Brenti,Gasharov,Neggers,Reiner-Welker,Stembridge-CE,Wagner1,Wagner2}. A labeled poset is {\em naturally labeled} 
if $x <_P y$ implies $\omega(x) < \omega(y)$. Neggers \cite{Neggers} made the conjecture for naturally labeled 
posets in 1978 and Stanley formulated the conjecture in its general form  in  1986. However, in 
\cite{Branden-CE}, we found a family of 
counterexamples to the Neggers-Stanley conjecture and subsequently Stembridge \cite{Stembridge-CE} found counterexamples that are naturally labeled thus disproving Neggers original conjecture. 

Although 
the Neggers-Stanley conjecture is refuted many questions regarding 
the  $(P,\omega)$-Eulerian polynomials remain open. A question which is still open is whether the 
coefficients of $W(P,\omega;t)$ always form a unimodal sequence. It is easy to see that real-rootedness 
implies unimodality. This weaker property was recently established by Reiner and Welker \cite{Reiner-Welker} 
 for a large and important class of posets, namely the class of naturally labeled and {\em graded} posets. A poset $P$ is graded if every saturated chain in $P$ has the same length. Prior to \cite{Reiner-Welker},  Gasharov \cite{Gasharov} proved unimodality for graded naturally labeled posets of rank at most $2$.   In \cite{Branden-Sign} 
 we proved unimodality for $(P,\omega)$-Eulerian polynomials of  labeled posets which we call {\em sign-graded posets}. The class of  sign-graded posets contains the class of naturally labeled graded posets. 
 


 If $(P,\omega)$ is a 
labeled poset we may associate signs to the edges of the Hasse-diagram, $E(P)$,  
of $P$ as follows. Let  
$\epsilon : E(P)  \rightarrow \{-1,1\}$ be defined by
$$
\epsilon(x,y)= \begin{cases} 
\ \ 1  \mbox{ if } \ \ \omega(x)<\omega(y), \\
-1 \mbox{ if } \ \ \omega(x)>\omega(y)   
\end{cases}
$$
It is not hard to prove that the $(P,\omega)$-Eulerian polynomial 
only depends on $\epsilon$, see \cite{Branden-Sign}. A labeled poset $(P,\omega)$ is 
sign-graded if for every maximal chain $x_1<x_2<\cdots < x_k$ in $P$, the 
sum of signs 
$$
\sum_{i=1}^k\epsilon(x_{i-1},x_i),
$$
is the same. Note that this definition extends the notion of 
graded posets since if $(P,\omega)$ is naturally labeled then 
all signs are equal to one and the above sum is just the length of 
the chain. The common value, $r$, of the above sum is called the {\em rank} 
of $(P,\omega)$. One may associate a (generalized) {\em rank function} 
$\rho : P \rightarrow \ZZ$ to a sign-graded poset by 
$$
\rho(x)=\sum_{i=1}^k\epsilon(x_{i-1},x_i),
$$
where  $x_1<x_2<\cdots < x_k=x$ is any saturated chain from 
a minimal element to $x$. In \cite{Branden-Sign} we prove the following theorem. 
\begin{theorem}[Br\"and\'en \cite{Branden-Sign}]\label{mainsign}
Let $(P,\omega)$ be a sign-graded poset of rank $r$ and let 
$d=p-r-1$. Then 
$$
W(P,\omega;t) = \sum_{i=0}^{\lfloor d/2 \rfloor} 
a_i(P,\omega)t^i(1+t)^{d-2i},
$$
where $a_i(P,\omega)$, $i=0,1,\ldots, \lfloor d/2 \rfloor$ are 
nonnegative integers. 
\end{theorem}

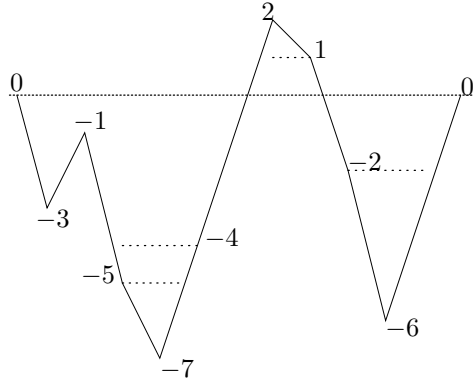
\begin{figure}\caption{\label{plusminus} The dotted lines indicates where 
the double ascents/descents are mapped.}
\bigskip
\setlength{\unitlength}{5mm}
\begin{picture}(13,12)
\put(0,0){
\path(0.2,8)(1,5)(2,7)(3,3)(4,1)(5,4)(7,10)(8,9)(9,6)(10,2)(12,8)
\put(0.02,8.1){$0$}\put(0.7,4.5){$-3$}\put(1.7,7.1){$-1$}\put(1.9,3){$-5$}\put(4,.5){$-7$}
\put(5.2,4){$-4$}\put(6.7,10){$2$}\put(8.1,9){$1$}\put(10,1.6){$-6$}\put(9,6){$-2$}\put(12,8){$0$}
\dottedline{0.2}(7,9)(8,9)
\dottedline{0.2}(9,6)(11,6)
\dottedline{0.2}(3,4)(5,4)
\dottedline{0.2}(3,3)(4.5,3)
\dottedline{0.1}(0,8)(12.3,8)
}
\end{picture}
\end{figure}

From the proof of Theorem~\ref{mainsign} in \cite{Branden-Sign} it is not evident  what the numbers $a_i(P,\omega)$ count. 
We will now give an alternative proof of Theorem~\ref{mainsign} by extending the MFS-action to $\JH(P,\omega)$. This will also give us an interpretation of the numbers $a_i(P,\omega)$,  
$i=0, \ldots, \lfloor d/2 \rfloor$. If both $(P,\omega)$ and $(P,\lambda)$ are 
sign-graded one can prove \cite[Corollary~2.4]{Branden-Sign} that 
up to a multiple of $t$ the corresponding  Eulerian polynomials 
are the same. Moreover, in \cite{Branden-Sign} we prove that  if $(P,\omega)$ is sign-graded then there exists a labeling $\mu$ of $P$ such that 
\begin{enumerate}
\item $(P,\mu)$ is sign-graded, 
\item the rank function of $(P,\mu)$ has values in $\{0,1\}$, 
\item all elements of rank $0$ have negative labels  and 
\item all elements of rank $1$ have positive labels
\end{enumerate}
Such a labeling will be called {\em canonical}. Hence it is no restriction in assuming that the sign-graded poset is labeled canonically. 
\begin{definition}
Let $(P,\omega)$ be sign-graded with $\omega$ canonical. For $x \in \omega(P)$ define a map 
$\psi_x : \JH(P,\omega) \rightarrow \JH(P,\omega)$ as follows. Let 
$\pi = a_1a_2 \cdots a_p \in \JH(P,\omega)$ and let 
$a_0=a_{p+1}=0$.
\begin{itemize}
\item If $x<0$ is a double descent let $\psi_x(\pi)$ be the permutation 
obtained by inserting $x$ between the first pair of letters 
$a_i, a_{i+1}$ to the right of $x$ such that $a_i < x < a_{i+1}$.
 
\item If $x<0$ is a double ascent let $\psi_x(\pi)$ be the permutation 
obtained by inserting $x$ between the first pair of letters 
$a_i, a_{i+1}$ to the left of $x$ such that $a_i > x > a_{i+1}$.

\item If $x>0$ is a double descent let $\psi_x(\pi)$ be the permutation 
obtained by inserting $x$ between the first pair of letters 
$a_i, a_{i+1}$ to the left of $x$ such that $a_i < x < a_{i+1}$

\item If $x>0$ is a double ascent  $\psi_x(\pi)$ be the permutation 
obtained by inserting $x$ between the first pair of letters 
$a_i, a_{i+1}$ to the right of $x$ such that $a_i > x > a_{i+1}$.    

\item If $x$ is a peak or a valley let $\psi_x(\pi)=\pi$.
\end{itemize}
See Fig.~\ref{plusminus}. 
\end{definition}
It is not immediate that this definition makes sense, i.e., that 
the resulting permutation represents a linear extension of $P$. 
Suppose that $x<0$ is a letter of $\pi \in \JH(P,\omega)$. Then $x$ is a letter of a maximal 
contigous subword $w$ of $\pi$ whose letters are all negative. By 
construction $\psi_x$ will not move $x$ outside of the 
word $w$. We claim that 
$$
\omega^{-1}(w)=\{y\in P: \omega(y) \mbox{ is a letter of } w\}
$$ 
is an anti-chain. Suppose that $y_1 <_P y_2$ are elements in  
$\omega^{-1}(w)$. Then, since $\rho(y_1)=\rho(y_2)=0$, there must 
be an element $z\in P$ such that $y_1 <_P z <_P y_2$, $\rho(z)=1$ and $\omega(z)>0$. 
This means that $\omega(z)$ is between $\omega(y_1)$ and $\omega(y_2)$ in $\pi$, so 
$\omega(z)$ is a letter of $w$ contrary to the assumption that 
all letters of $w$ are negative. Since $\omega^{-1}(w)$ is an anti-chain and since $\psi_x$ does not move $x$ outside $\omega^{-1}(w)$ we have that $\psi_x(\pi) \in \JH(P,\omega)$. The case 
$x>0$ is analogous. 

We may now define a $\ZZ_2^P$-action on $\JH(P,\omega)$ by  
$$
\psi_S(\pi)= \prod_{x \in S}\psi_{\omega(x)}(\pi), \ \ \ \ S \subseteq P.
$$  
Let 
$\hat{\pi}$ be the unique permutation in $\Orb(\pi)$ such 
that $0\hat{\pi}0$ has no double descents.  

\begin{theorem}\label{WP}
Let $(P,\omega)$ be a sign-graded poset of rank $r$ 
where $\omega$ is canonical and let  
$\pi \in \JH(P,\omega)$. Then  
$$
\sum_{\sigma \in \Orb(\pi)} t^{\des(\sigma)} =
t^{\des(\hat{\pi})} (1+t)^{p-r-1-2\des(\hat{\pi})}.
$$
Moreover, if $r=0$ then $\peak(\cdot)$ is invariant under the $\ZZ_2^P$-action and 
$\peak(\pi)=\des(\hat{\pi})$ for all $\pi \in \JH(P,\omega)$.  
\end{theorem}
\begin{proof}
If $x$ is a double ascent in $0\pi0$ then $\des(\psi_x(\pi))= \des(\pi)+1$. It follows that 
$$
\sum_{\sigma \in \Orb(\pi)} t^{\des(\sigma)} = t^{\des(\hat{\pi})}(1+t)^a 
$$ 
where $a$ is the number of double ascents in $\pi$. Suppose $r=0$. Deleting all 
double ascents in $\hat{\pi}$ results in an alternating permutation 
$$
0>a_1<a_2>a_3< \cdots >a_{p-a}<0,
$$
with the same number of peaks/descents as $\pi$. Hence $p-a=2\peak(\pi)+1$. 

If $r=1$, deleting all 
double ascents in $\hat{\pi}$ results in an alternating permutation 
$$
0>a_1<a_2>a_3< \cdots  <a_{p-a}>0,
$$
with the same number of descents. Hence $p-a-2=2\des(\hat{\pi})$. 
\end{proof}
Stembridge \cite{Stembridge-peak} developed a theory of "enriched $P$-partitions" in which the distribution of {\em peaks} in $\JH(P,\omega)$ and the polynomial, viz., 
$$
\WP(P,\omega;t) = \sum_{\pi \in \JH(P,\omega)}t^{\peak(\pi)},
$$
play a significant role. 
For a canonically labeled poset $(P,\omega)$ let $(\hat{P},\hat{\omega})$ be any canonically labeled 
poset such that $\hat{P}$ is obtained from $P$ by adjoining a greatest element. 
\begin{theorem}
Let $(P,\omega)$ be a canonically labeled sign-graded poset of rank $r$. If $r=0$ then 
$$
W(P,\omega;t) = 2^{-p+1}(1+t)^{p-1}\WP\big(P,\omega; 4t(1+t)^{-2}\big).
$$ 
Equivalently, 
$$
a_i(P,\omega)= 
2^{-p+1+2i} |\{ \pi \in \JH(P,\omega): \peak(\pi)=i\}|. 
$$

\noindent
If $r=1$ then 
$$
W(P,\omega;t) = 2^{-p}t^{-1}(1+t)^{p}\WP\big(\hat{P},\hat{\omega}; 4t(1+t)^{-2}\big).
$$ 
Equivalently, 
$$
a_i(P,\omega)= 2^{-p+2+2i}|\{ \pi \in \JH(\hat{P},\hat{\omega}): \peak(\pi)=i+1\}|. 
$$
\end{theorem}
\begin{proof}
Note that $W(\hat{P},\hat{\omega};t)=t^{-r}W(P,\omega;t)$, so we may assume that $r=0$. By 
Theorem~\ref{WP} the proof follows just as  the proof of Corollary~\ref{corre}. 
\end{proof}

\section{Vertices Of Even Height}\label{veh}
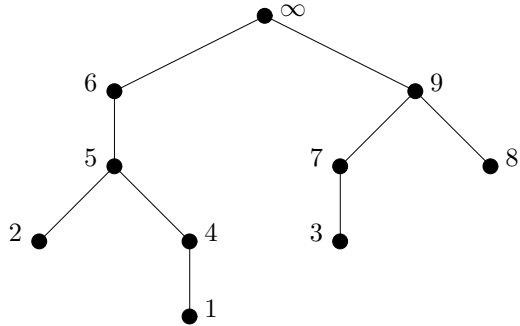
\begin{figure}\caption{\label{dectree} The decreasing unordered tree corresponding to 
$652419738$}
\setlength{\unitlength}{20mm}
\bigskip\bigskip
\newcommand\p{\circle*{0.1}}
\begin{picture}(3,1.5)
\put(0,0){
\path(2,1.5)(1,1)\put(1,1){$\p$}\put(0.8,1){$6$}
\put(2,1.5){$\p$}\put(2.1,1.5){$\infty$}
\path(2,1.5)(3,1)\put(3,1){$\p$}\put(3.1,1){$9$}
\path(3,1)(3.5,0.5)\put(3.5,0.5){$\p$}\put(3.6,0.5){$8$}
\path(3,1)(2.5,0.5)\put(2.5,0.5){$\p$}\put(2.3,0.5){$7$}
\path(2.5,0.5)(2.5,0)\put(2.5,0){$\p$}\put(2.3,0){$3$}
\path(1,1)(1,0.5)\put(1,0.5){$\p$}\put(0.8,0.5){$5$}
\path(1,0.5)(0.5,0)\put(0.5,0){$\p$}\put(0.3,0){$2$}
\path(1,0.5)(1.5,0)\put(1.5,0){$\p$}\put(1.6,0){$4$}
\path(1.5,0)(1.5,-0.5)\put(1.5,-0.5){$\p$}\put(1.6,-0.5){$1$}
}
\end{picture}
\bigskip\bigskip\bigskip\bigskip 
\end{figure}

To any permutation $w$ of a finite subset of $\{1,2,\ldots\}$ we may associate a decreasing 
unordered tree as follows. Let 
$\infty$ be a symbol which is greater than every letter in $w$. 
If $w$ is empty then $T(w;\infty)$ is the tree with a single vertex 
labeled $\infty$. Otherwise write $w$ as 
$w=m_1w_1m_2w_2 \cdots m_kw_k$ where $m_i$ are the left-to-right maxima 
of $w$. Then $T(w;\infty)$ is the labeled tree with $T(w_i;m_i)$ as 
subtrees of the root, see Fig.~\ref{dectree}.    
Let $\veh(\pi)$ be the 
number of (non-root) vertices of even height in $T(\pi;\infty)$. As Fig.~\ref{dectree} suggests 
$$
\veh(652419738)=|\{1,5,7,8\}|=4. 
$$ 
We will here show that $\veh$ and 
$\des$ have the same distribution on any subset of $\sym_n$ invariant under the MFS-action.  For 
$\pi \in \sym_n$ and $x \in [n]$ let $r_\pi(x)$ be the number of right edges in the path from the root to $x$ in the decreasing binary tree associated with $\pi$. It is plain to see that $r_\pi(x)+1$ is equal to the height of $x$ as a vertex  $T(\pi;\infty)$. Let 
$\Odd(\pi)$ the set of all $x \in [n]$ for which $r_\pi(x)$ is odd. Hence $\Odd(\pi)$ is the set of vertices of even height in $T(\pi;\infty)$.  Also, let $\Redge(\pi)$ be the set of vertices in the decreasing binary tree which are ends of right edges. Clearly, $\des(\pi)=|\Redge(\pi)|$.  Define  $\Psi, \Phi : \sym_n \rightarrow \sym_n$ by 
\begin{eqnarray*}
\Psi(\pi) &=& \prod_{x} \varphi_x(\pi) \quad (x \in \Odd(\pi));\\
\Phi(\pi) &=&  \prod_{x} \varphi_x(\pi) \quad (x \in \Redge(\pi)).
\end{eqnarray*}

\begin{theorem}\label{psiphi}
The transformations $\Psi$ and $\Phi$ are inverses of each other. Moreover, if $\pi \in \sym_n$ then 
\begin{eqnarray*}
\Odd(\pi) &=& \Redge(\pi') \quad \mbox{and } \\
\Redge(\pi) &=& \Odd(\pi''), 
\end{eqnarray*}
where $\pi'=\Psi(\pi)$ and $\pi'' = \Phi(\pi)$.
\end{theorem}
\begin{proof}
Note that it is enough to prove the first equality since then 
\begin{eqnarray*}
\Phi (\Psi (\pi)) &=& \prod_{x}\varphi_x \prod_{y} \varphi_y (\pi)     \quad    \big(x \in \Redge(\pi'), y \in \Odd(\pi)\big)\\  
                        &=&  \prod_{x}\varphi_x \prod_{y} \varphi_y(\pi)   \quad  \big(x,y \in \Odd(\pi)\big) \\
 &=&\pi,
 \end{eqnarray*}
 because the involutions $\varphi_x$ commute.

Let $T$ and $T'$ be the decreasing binary trees corresponding to $\pi$ and $\pi'$, respectively. 
Suppose that $x \in \Odd(\pi)$ and let $y$ be its father in $T$. If there is a right-edge between $x$ and $y$ in $T$ then 
$y \notin \Odd(\pi)$, which means that  there will also be a right-edge between $x$ and $y$ in $T'$, so that $x \in \Redge(\pi')$. If there is a left-edge between $x$ and $y$ in $T$ then also $y \in \Odd(\pi)$ so 
that $\Psi$ will turn this edge to a right-edge and hence $x \in \Redge(\pi')$. 

The fact that  $x \notin \Odd(\pi)$ implies $x \notin \Redge(\pi)$ follows similarly. 
\end{proof}
Note that $\Psi$ and $\Phi$ restricts to bijections on  all subsets of $\sym_n$ invariant under 
the "proper" Foata-Strehl action, but not on subsets invariant under the modified Foata-Strehl action. 
Define a transformation $\Psi' : \sym_n \rightarrow \sym_n$ by 
$$
\Psi'(\pi) = \prod_{x} \varphi'_x(\pi) \quad (x \in \Odd(\pi)).
$$
\begin{theorem}\label{psi'}
Let $T \subseteq \sym_n$ be invariant under the modified Foata-Strehl action. Then 
$\Psi' : T \rightarrow T$ is a bijection and 
$$
\veh(\pi) = \des( \Psi'(\pi) ),  \quad  \pi \in T.
$$
\end{theorem}
\begin{proof}
Since the involutions $\varphi_x$ commute we may write $\Psi'$ as $\Psi' = F \circ \Psi$ where 
$F$ is defined by 
$$
F(\pi) = \prod_x \varphi_x(\pi) \quad (x \in \Redge(\pi), c(x)=2),
$$ 
and where $c(x)$ is the number children of $x$ in the decreasing binary tree of $\pi$. Clearly, 
$\des(\pi) = \des(F(\pi))$ so it remains to prove that $\Psi'$ is a bijection. 

Let $f$ be defined as in Remark~\ref{comb} and let $\pi \in \sym_n$. Then since the involutions 
$\varphi_x$ commute we have 
\begin{eqnarray*}
f(\Psi'(\pi))&=&  \prod_{y} \varphi'_y \prod_x \varphi'_x(\pi) \quad (y \in [n], x \in \Odd(\pi)) \\
&=& \prod_x\varphi'_x(\pi) \quad (x \notin \Odd(\pi)).
\end{eqnarray*}

It follows that $\Psi'$ can be defined 
recursively on the set of permutations of any finite subset of $\{1,2,\ldots\}$ as follows. The 
empty word is mapped by $\Psi'$ to itself, and if $w = LnR$ where $n$ is 
the greatest element of $w$ and $L$ and $R$ are the words to the left and right of $n$ respectively then 
$$
\Psi'(w) = \Psi'(L)nf(\Psi'(R)), 
$$
where  $f$ is as in Remark~\ref{comb}. From this recursive definition it is plain to see that $\Psi'$ is bijective.  
\end{proof}

\begin{corollary}
Let $n,r \in \NN$. Then $\veh$ and $\des$ have the same distribution over $\sym_n^r$.
\end{corollary}

To every unordered decreasing tree $T(\pi; \infty)$ corresponding to a permutation 
$\pi \in \sym_n(231)$ there is a unique ordered unlabeled tree obtained 
by ordering the children of a vertex decreasingly from left to right and
dropping the labels. Recall that a {\em Dyck-path} of length 
$2n$ is a lattice path in $\NN^2$ starting at the origon and ending 
at $(2n,0)$, using steps $u = (1,1)$ and $d=(1,-1)$, and never 
going below the $x$-axis. If we traverse the ordered tree in pre-order 
and write a $u$ every time we go down an edge and write a $d$ every 
time we go up an edge we obtain a Dyck path. This 
describes a bijection between the set of Dyck path of length $2n$ and 
the set of ordered trees with $n+1$ vertices (and by the above also between the set of Dyck path of length $2n$ and $\sym_n(231)$). Note that a  vertex of even 
height translates into an up-step of even height in the Dyck path,  
and a descent translates into a double up-step $uu$ in the path. We  have 
thus recovered the following classical result of Kreweras \cite{Kreweras}. 
\begin{corollary}
The statistics "up-steps at even height" and "double up-steps" have the same  
distribution over the set of Dyck paths of a given length. 
\end{corollary}  
When restricted to $\sym_n(231)$ one may express $\veh$ as the following alternating sum of 
permutation patterns \cite{Branden-C-S} 
$$
\veh(\pi)= d_1(\pi)-2d_2(\pi)+4d_3(\pi)- \cdots + (-2)^{n-2}d_{n-1}(\pi),
$$
where $d_i(\pi)$ is the number of decreasing subsequences of length $i+1$ in $\pi$.  

\section{A Mahonian Partner for Vertices of Even Height}\label{mahonian}
\begin{figure}\caption{\label{tree} The increasing unordered tree corresponding to 
$586317492$}
\setlength{\unitlength}{20mm}
\bigskip\bigskip
\newcommand\p{\circle*{0.1}}
\begin{picture}(3,1.5)
\put(0,0){
\path(2,1.5)(1,1)\put(1,1){$\p$}\put(0.8,1){$1$}
\put(2,1.5){$\p$}\put(2.1,1.5){$0$}
\path(2,1.5)(3,1)\put(3,1){$\p$}\put(3.1,1){$2$}
\path(3,1)(3.5,0.5)\put(3.5,0.5){$\p$}\put(3.6,0.5){$9$}
\path(3,1)(2.5,0.5)\put(2.5,0.5){$\p$}\put(2.3,0.5){$4$}
\path(2.5,0.5)(2.5,0)\put(2.5,0){$\p$}\put(2.3,0){$7$}
\path(1,1)(1,0.5)\put(1,0.5){$\p$}\put(0.8,0.5){$3$}
\path(1,0.5)(0.5,0)\put(0.5,0){$\p$}\put(0.3,0){$5$}
\path(1,0.5)(1.5,0)\put(1.5,0){$\p$}\put(1.6,0){$6$}
\path(1.5,0)(1.5,-0.5)\put(1.5,-0.5){$\p$}\put(1.6,-0.5){$8$}
}
\end{picture}
\bigskip\bigskip\bigskip\bigskip
\end{figure}
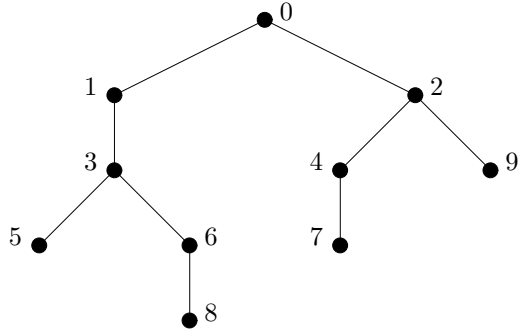

Recall that the {\em descent set} of a permutation $\pi=a_1a_2\cdots a_n$ is defined by 
$\Des(\pi)=\{i \in [n-1]: a_i> a_{i+1}\}$ and that the {\em major index} of $\pi$ as 
$$
\MAJ(\pi) = \sum_{i \in \Des(\pi)}i. 
$$  
A statistic $B : \sym_n \rightarrow \NN$ is said to be {\em Mahonian} if it has the same distribution as 
$\MAJ$ on $\sym_n$, i.e., 
$$
\sum_{\pi \in \sym_n}q^{B(\pi)}= [n]_q[n-1]_q \cdots [1]_q, 
$$
where $[k]_q=1+q+\cdots+q^{k-1}$. A bi-statistic $(A,B)$ is {\em Euler-Mahonian} if it has the same distribution as 
$(\des,\MAJ)$ on $\sym_n$. 
We will now redefine the statistic $\veh$ so that 
we can define a Mahonian partner for it. To every permutation $\pi=a_1a_2\cdots a_n \in \sym_n$ we 
associate an increasing unordered tree,  $T'(\pi)$, as follows.   
If $b$ is a right-to-left minimum of $\pi$ then $b$ is a successor of the root, which is labeled $0$. 
Otherwise $b$ is the successor of the 
leftmost element $a$ to the right of $b$ which is smaller than 
$b$, see Fig.~\ref{tree}. 

Let $n \in \NN$. We (re-)define the statistic {\em vertices of even height}, $\vej: \sym_n \rightarrow \NN$, by letting $\vej(\pi)$ be the number of (non-root) vertices in $T'(\pi)$ of even
height. Thus $\vej (586317492) = |\{3,4,8,9\}|=4$. 
We define the {\em even vertex set}, $\EV(\pi)$, as the set of indices 
$1 \leq i \leq n$ 
such that the vertex $a_i$ is of even height.

The {\em complement}, $\pi^c$, of a permutation $\pi=a_1a_2\cdots a_n$ is the 
permutation $b_1b_2 \cdots b_n$ on the same letters as $\pi$ such that $a_i < a_j$ if and only if $b_i>b_j$ for all
$1 \leq i < j \leq n$. We define a 
transformation $\Theta$ on permutations of any finite subset of $\{1,2,\ldots\}$ 
recursively as follows. The empty permutation is mapped  
onto itself and if $\pi$ is the concatenation $\sigma m \tau$ where 
$m$ is the smallest letter in $\pi$,  then $\Theta(\pi) = \Theta(\sigma^c) m \Theta(\tau)$. It is clear 
that $\Theta$ restricted to the symmetric group is a bijection. 
\begin{eqnarray*}
\Theta(586317492)&=&\Theta(6358)1\Theta(7492)\\
               &=&63\Theta(58)1\Theta(794)2\\
               &=&635819742.
\end{eqnarray*} 
If 
$S \subset \ZZ$ and $k \in \ZZ$ let $S+k := \{ s+k : s \in S\}$. 
\begin{theorem}\label{evt}
Let $n \in \NN$. For all permutations $\pi \in \sym_n$ we have
$$\EV(\Theta(\pi)) = \Des(\pi).$$
\end{theorem}

\begin{proof}
The proof is by induction over the length $n=|\pi|$ of $\pi$. The case $n=0$ is 
clear. Suppose that $n>0$. Then we can write $\pi \in \sym_n$ as the 
concatenation $\sigma 1 \tau$. 
Let $k = |\sigma|$. If $1 \leq i \leq k$ then clearly $i \in \EV(\sigma)$ if and 
only if $i \notin \EV(\pi)$.
Hence 
\begin{eqnarray*}
\EV(\pi) &=& \Big([k]\setminus \EV(\sigma)\Big)\cup \Big( \EV(\tau)+k+1\Big) \ \ \ \ \mbox{ and }\\
\Des(\pi)
         &=& \Big([k]\setminus \Des(\sigma^c)\Big)\cup \Big(\Des(\tau)+k+1\Big),
\end{eqnarray*}
since $[n] \setminus \Des(\pi) = \{n\}\cup \Des(\pi^c)$ for all $\pi$ of length $n$. Using induction we get 
\begin{eqnarray*}
\EV(\Theta(\pi)) &=& \EV\Big(\Theta(\sigma^c)1\Theta(\tau)\Big) \\
               &=& \Big([k]\setminus \EV(\Theta(\sigma^c))\Big)\cup\Big(\EV(\Theta(\tau))+k+1\Big)\\
             &=& \Big([k] \setminus \Des(\sigma^c)\Big)\cup\Big(\Des(\tau)+k+1\Big)\\
             &=& \Des(\pi).            
\end{eqnarray*}
\end{proof}
It is desirable to find a bijection which is not defined recursively and which proves Theorem~\ref{evt}.  

We may now define a Mahonian partner for $\vej$. 
The statistic sum of indices of vertices even height,  
$\SIVEH : \sym_n \rightarrow \NN$, 
is defined by
$$ \SIVEH(\pi) = \sum_{i \in \EV(\pi)}i .$$
\begin{corollary}\label{mahonian}
For all $n \in \NN$ the bistatistic $(\vej, \SIVEH)$ is Euler-Mahonian on $\sym_n$.
\end{corollary}

\section{Gal's Conjecture on $\gamma$-Polynomials}\label{gamma}
Recall that the $h$-polynomial of a simplicial complex $\Delta$ of dimension $d-1$ is the polynomial $h_\Delta(t)= h_0(\Delta)+h_1(\Delta)t+\cdots+h_d(\Delta)t^d$ defined by the polynomial identity 
$$
\sum_{i=0}^dh_i(\Delta)t^i(1+t)^{d-i}= \sum_{i=0}^df_{i-1}(\Delta)t^i, 
$$
where $f_i(\Delta)$, $-1\leq i \leq d-1$ is the number of faces of $\Delta$ of dimension $i$. If $\Delta$ is a simplicial homology sphere then the {\em Cohen-Macaulay property} and the {\em Dehn-Sommerville equations} imply that $\{h_i(\Delta)\}_{i=0}^d$ is nonnegative and symmetric.  Hence one may define the $\gamma$-polynomial of $\Delta$,  
$\gamma_\Delta (t)= \sum_{i=0}^{\lfloor d/2 \rfloor} \gamma_i(\Delta)t^i$, by 
$$
h_\Delta(t)=  \sum_{i=0}^{\lfloor d/2 \rfloor} \gamma_i(\Delta)t^i(1+t)^{d-2i}.
$$
A simplicial complex  
$\Delta$ is {\em flag} if  the minimal non-faces of $\Delta$ have cardinality two. The following conjecture 
generalizes the Charney-Davis conjecture \cite{Charney}. 
\begin{conjecture}[Gal \cite{Gal}]\label{gal}
If $\Delta$ is a flag simplicial homology sphere of dimension $d-1$, then 
$$
\gamma_i(\Delta) \geq 0, \ \ \ \   0 \leq i \leq \lfloor d/2 \rfloor. 
$$
\end{conjecture}
It is desirable to find a combinatorial, geometrical or ring-theoretical description of the numbers 
$\gamma_i(\Delta)$. In \cite{Reiner-Welker} Reiner and Welker associated to any graded naturally labeled poset $(P,\omega)$ a simplicial polytopal sphere, $\Delta_{eq}(P)$, whose $h$-polynomial is the $(P,\omega)$-Eulerian polynomial. Hence, Theorem~\ref{WP} gives a combinatorial description of the $\gamma$-polynomial of $\Delta_{eq}(P)$ and verifies Conjecture~\ref{gal} for $\Delta_{eq}(P)$.  

In \cite{Postnikov-RW} Postnikov, Reiner and Williams extended the MFS-action to give a combinatorial interpretation of the $\gamma$-polynomials of tree-associahedra which confirms Conjecture~\ref{gal} in this case. Also, Chow \cite{Chow} has  given  a combinatorial interpretation of the $\gamma$-polynomials of the Coxeter complexes of type $B$ and $D$ and confirming Conjecture~\ref{gal} for these complexes. 

\section{Further Directions and Open Problems}\label{OP}
Let $\mathcal{I}_n$ be the set of involutions in $\sym_n$ and let 
$$
I_n(t)=\sum_{\pi \in \mathcal{I}_n}t^{\des(\pi)}=\sum_{k=0}^{n-1}I_{n,k}t^k. 
$$
Brenti has conjectured that the sequence $\{I_{n,k}\}_{k=0}^{n-1}$ has no internal zeros and is log-concave, i.e., 
$$
I_{n,k}^2 \geq I_{n,k+1}I_{n,k-1},  \ \ 1 \leq k \leq n-2,  
$$
see \cite{Dukes} where progress on this conjecture was made. Motivated by Brenti's conjecture Guo and Zeng \cite{Guo} proved the weaker statement that $\{I_{n,k}\}_{k=0}^{n-1}$ is unimodal. Also, Strehl \cite{Strehl} proved symmetry for $\{I_{n,k}\}_{k=0}^{n-1}$ and 
the following conjecture was made in \cite{Guo}.   
\begin{conjecture}[Guo-Zeng \cite{Guo}]\label{GZ}
Let $n \in \NN$. Then 
$$
I_n(t) = \sum_{i=0}^{\lfloor (n-1)/2 \rfloor} a_{n,i}t^i(1+t)^{n-1-2i}, 
$$
where $a_{n,i} \in \NN$ for $0 \leq i \leq  \lfloor (n-1)/2 \rfloor$. 
\end{conjecture} 
Gessel \cite{Gessel} has conjectured a fascinating property of the joint distribution of descents and inverse descents. 
\begin{conjecture}[Gessel \cite{Gessel}]\label{gesselc}
Let $\tau \in \sym_n$. Then 
\begin{equation}\label{gesseleq}
\sum_{\pi \in \sym_n}s^{\des(\pi)}t^{\des(\pi^{-1}\tau)}= 
\sum_{k,j}c_n(\tau;k,j)(s+t)^k(st)^j(1+st)^{n-k-1-2j},
\end{equation}
where $c_n(\tau;k,j) \in \NN$ for all $k,j \in \NN$.
\end{conjecture}
Symmetry properties imply that an expansion such as \eqref{gesseleq} with $c_n(\tau;k,j) \in \ZZ$, $k,j \in \NN$ exists. Moreover, $c_n(\tau;k,j)$ only depends on the number of descents of $\tau$. In light of  Conjectures~\ref{GZ} and \ref{gesselc} there might be another $\ZZ_2^n$-action on permutations which also behaves well with respect to the inverse permutation.  

Recall the definition of $A_n(p,q,t)$ of Section~\ref{pq}. The first nontrivial examples are
\begin{eqnarray*}
A_3(p,q,t) &=&(1+t)^2 + (p+q)t \\
A_4(p,q,t) &=&(1+t)^3 + (p+q)(p+q+2)t(1+t) \\
A_5(p,q,t) &=&(1+t)^4 + (p+q)\left((p+q)^2+2(p+q)+3)\right)t(1+t)^2 + \\
    && (p+q)^2(p^2 + pq + q^2 +1)t^2.
\end{eqnarray*}

\begin{conjecture}
Let $b_{n,i}(q)$ be defined by \eqref{bni}. Then  
$(p+q)^i \mid b_{n,i}(p,q)$ for all 
$0 \leq i \leq \lfloor (n-1)/2 \rfloor$. 
\end{conjecture}

%
%
\bigskip

\begin{center}
{  A{\Small CKNOWLEDGEMENTS}} 
\end{center}

The author would like to thank the anonymous referee for helpful advice on the presentation of this paper, for pointing out the correct credit for the actions and for showing how the bijection $\Psi'$ in 
Theorem~\ref{psi'} can be described non-recursively, as presented in Section~\ref{veh}.

\end{document}